\documentclass[12pt,draft]{amsart}

\usepackage{xypic}
\usepackage{amssymb}

\newtheorem{lemma}{Lemma}[section]

\newtheorem{proposition}[lemma]{Proposition}
\newtheorem{theorem}[lemma]{Theorem}
\newtheorem{conjecture}[lemma]{Conjecture}
\newtheorem{corollary}[lemma]{Corollary}

\theoremstyle{definition}
\newtheorem{defn}[lemma]{Definition}

\theoremstyle{remark}
\newtheorem{rmk}[lemma]{Remark}
\newtheorem{example}[lemma]{Example}
\newtheorem{remark}[lemma]{Remark}

\newcommand{\Z}{\ensuremath{{\mathbb Z}}}
\newcommand{\Q}{\ensuremath{\mathbb Q}}

\renewcommand{\div}{{\rm Div}}

\newcommand{\isoarrow}{\tilde{\longrightarrow}}

\newcommand{\spec}{{\rm Spec}}

\newcommand{\A}{{\mathbb A}}
\newcommand{\Gm}{{\mathbb G}_m}
\newcommand{\R}{{\mathbb R}}
\newcommand{\cc}{{\mathbb C}}
\newcommand{\uU}{{\mathfrak u}}
\newcommand{\qq}{{\mathbb Q}}
\newcommand{\ffq}{{\mathbb F_q}}
\newcommand{\qql}{{\mathbb Q_{\ell}}}
\newcommand{\qqlbar}{{\overline{\mathbb Q}_\ell}}
\newcommand{\pp}{{\mathbb P}}
\newcommand{\G}{{\mathbb G}}
\newcommand{\zz}{{\mathbb Z}}
\newcommand{\Char}{\mathop{\rm char}}
\newcommand{\tr}{\mathop{\rm tr}}
\newcommand{\ol}{\overline}
\newcommand{\Zeta}{Z}
\renewcommand{\O}{{\mathcal O}}
\newcommand{\vol}{{\mathop{\rm vol}}}

\newcommand{\sym}[2][n]{{#2}^{(#1)}}
\newcommand{\kk}{{\mathfrak K}}

\renewcommand{\div}{{\rm Div}}
\newcommand{\rx}{{\rm X}}
\newcommand{\dom}{{{\rm dom}}}
\newcommand{\Aut}{\mathop{\rm Aut}}
\newcommand{\Hom}{\mathop{\rm Hom}}
\newcommand{\codim}{\mathop{\rm codim}\nolimits}

\newcommand{\sX}{{\mathfrak X}}
\newcommand{\sY}{{\mathfrak Y}}
\newcommand{\sZ}{{\mathfrak Z}}

\newcommand{\bun}[2]{{\mathfrak{Bun}}^{#1}_{#2}}

\newcommand{\Var}{\mathop{V\!ar}\nolimits}
\renewcommand{\L}{{\mathbb L}}
\newcommand{\kvar}{K_0(\Var_k)}
\newcommand{\fqvar}{K_0(\Var_{\ffq})}
\newcommand{\kvarhat}{\widehat{K}_0(\Var_k)}
\newcommand{\cvarhat}{\widehat{K}_0(\Var_{\cc})}
\newcommand{\fqvarhat}{\widehat{K}_0(\Var_{\ffq})}

\newcommand{\kvarhatG}{\widehat{K}_0^G(\Var_k)}

\newcommand{\gln}[1][n]{{\rm GL}_{#1}}
\newcommand{\sln}[1][n]{{\rm SL}_{#1}}

\newcommand{\sm}{{\rm Sm}}
\newcommand{\cor}{{\rm SmCor}}
\newcommand{\dmgm}[1][\Q]{{\rm DM}^{\rm eff}_{\rm gm}(k,#1)}
\newcommand{\dmm}[1][\Q]{{\rm DM}^{\rm eff}_{-}(k,#1)}
\newcommand{\cs}{{\rm C}_*}
\newcommand{\eumc}{{\chi}^c_{\rm Mot}}

\begin{document}
\sloppy
\title{On the Motive of the Stack of Bundles}

\author{Kai Behrend}
\address{University of British Columbia}
\email{behrend@math.ubc.ca}

\author{Ajneet Dhillon}
\address{University of Western Ontario}
\email{adhill3@uwo.ca}

\begin{abstract}
Let $G$ be a split connected semisimple group over a field. 
We give a conjectural formula for the motive of the stack of
$G$-bundles over a curve $C$, in terms of special values of the
motivic zeta function of $C$.  The formula is true if  $C=\pp^1$ or
$G=\sln$. If $k=\cc$, upon applying the Poincar\'e or Serre
characteristic, the 
formula reduces to results of Teleman and Atiyah-Bott on  the gauge
group. If $k=\ffq$, upon applying the counting 
measure, it reduces to the fact that the Tamagawa number of
$G$ over the function field of $C$ is $|\pi_1(G)|$.
\end{abstract}

\maketitle

\section{Introduction}

We work over a ground field $k$.
For a variety $Y$ we write $\mu(Y)$ for its class in the
$K$-ring of varieties, $\kvar$.

As any principal $\gln$-bundle (or $\gln$-torsor) $P\rightarrow X$
($X$ a variety) is
locally trivial in the Zariski topology, we
have the formula $\mu(P) = \mu(X)\mu(\gln)$.
We will use this fact to define $\mu(\sX)\in\kvarhat$ whenever $\sX$
is an algebraic stack stratified by global quotients.  Here $\kvarhat$
is the dimensional completion of $\kvar[\frac{1}{\L}]$, in which
$\mu(\gln)$ is invertible. In fact, if $\sX\cong[X/\gln]$ is a global
quotient, we define $\mu(\sX)=\frac{\mu(X)}{\mu(\gln)}$, and
generalize from there. Note that every Deligne-Mumford stack of
finite type  is stratified by global quotients.

We will also introduce a variation on $\kvarhat$, namely the
modified ring $\kvarhatG$ obtained by imposing the extra relations
(the `torsor relations')
\[
\mu(P) = \mu(X)\mu(G),
\]
whenever $P\rightarrow X$ is a $G$-torsor and $G$ is a fixed connected split
linear algebraic group. One can show that all the usual
characteristics factor through this ring. In the appendix, the second
author will show that there is a ring homomorphism
\[
\kvarhatG\rightarrow {K}_0(\dmgm),
\]
where $\dmgm$ is the $\Q$-linearization of Voevodsky's
triangulated category of effective geometrical motives
and $\Char(k)=0$.

Throughout this paper $C$ is a smooth projective geometrically
connected curve over $k$.  We fix also a split semisimple connected
algebraic group $G$ over $k$.  Let $\bun{}{G,C}$ denote the moduli
stack of $G$-torsors on $C$. The stack $\bun{}{G,C}$ is stratified by
global quotients, and even though it is not of finite type, its motive
still converges in $\kvarhat$, because the dimensions of the boundary
strata (where the bundle becomes more and more unstable) tend to
$-\infty$.

The purpose of this paper is to propose a conjectural formula for the
motive of $\bun{}{G,C}$ in $\kvarhatG$. Our formula expresses
$\mu(\bun{}{G,C})$ in 
terms of special values of the motivic zeta function of $C$.
For simply connected $G$, the formula reads:
$$\mu(\bun{}{G,C})= \L^{(g-1)\dim G}
\prod_{i=1}^r \Zeta(C,\L^{-d_i})\,,
$$
where the $d_i$ are the
numbers one higher than the exponents of $G$.

If $k$ is a finite field, we can apply the counting measure to this formula.
We obtain a statement equivalent to the celebrated conjecture of Weil, 
to the effect that the
Tamagawa number of $G$ (as a group over the function field of $C$)
is equal to~1. Of course, Weil's
conjecture  is much more general, as it  applies to arbitrary 
semisimple simply connected groups over any  
global field.

The proof of the Tamagawa number conjecture in the case of a 
split group induced from the ground field 
was completed by
Harder \cite{harder:74} 
by studying residues of Eisenstein series and using an idea
of Langlands. Motivic Eisenstein series have been defined
in \cite{kapranov:00} so it is natural to ask if there
is a proof of our conjecture along similar lines.

We consider our conjecture to be  a motivic version of Weil's 
Tamagawa number conjecture.
Thus we are lead to consider
$$\tau(G)=\L^{(1-g)\dim G}\mu(\bun{}{G,C})\prod_{i=1}^r\Zeta(C,\L^{-d_i})^{-1}
\quad \in \quad\kvarhatG$$
as the {\em motivic Tamagawa number }of $G$. We hope to find an
interpretation of $\tau(G)$ as a measure in a global motivic
integration theory, to be developed in the future.

We provide four pieces of evidence for our  conjecture:

In Section~\ref{gauge}, we prove that if $k=\cc$ and we apply the
Poincar\'e characteristic to our conjecture, the simply connected case
is true. It follows from  results on the
Poincar\'e series of the gauge 
group of $G$ and  the purity of the Hodge structure of
$\bun{}{G,C}$ due to Teleman~\cite{teleman:98}.

In Section~\ref{auto}, we verify that if $k=\ffq$, and we apply the
counting measure to our conjecture it reduces to theorems of Harder
and Ono that assert that the
Tamagawa number of $G$ is the cardinality of the fundamental group of
$G$. 

In Section~\ref{sln}, we prove our conjecture for $G=\sln$ using
the construction of matrix divisors in \cite{bifet:94}.

Finally, in Section~\ref{p1}, we prove our conjecture for $C=\pp^1$,
using  the explicit classification of $G$-torsors due
to Grothendieck and Harder.

\section{The Motive of an Algebraic Stack}

\subsection{Dimensional completion of the $K$-ring of varieties}

Let $k$ be a field.    The underlying abelian group of the 
ring $\kvar$ is  generated by the symbols
$\mu(X)$, where $X$ is the isomorphism class of a variety over
$k$, subject to all relations
$$
\mu(X) = \mu(X\setminus Z) + \mu(Z)\quad \text{ if $Z$ is closed
in $X$}\,.
$$
We   call $\mu(X)$ the {\em motive }of $X$.

Cartesian product of varieties induces a ring structure on
$\kvar$.  Thus $\kvar$ becomes a commutative ring with unit. Let
$\L$ denote the class of the affine line in $\kvar$.

The ring $\kvarhat$ is obtained by taking the dimensional
completion of $\kvar$. Explicitly, define $F^m(\kvar_\L)$ to
be the abelian subgroup of $$\kvar_\L=\kvar[\textstyle{\frac{1}{\L}}]$$
generated by symbols of the form
\[
\frac{\mu(X)}{\L^n}
\]
where $\dim X - n \le -m$. This is a ring filtration and
$\kvarhat$ is obtained by completing $\kvar_\L$ with respect to
this filtration.

Note that $\L^n-1$ is invertible in $\kvarhat$ as
\begin{eqnarray*}
  \frac{1}{\L^n - 1} &=& \frac{1}{\L^{n}}(\frac{1}{1-\frac{1}{\L^n}}) \\
                     &=& \L^{-n}(1 + \L^{-n} + \ldots ).
\end{eqnarray*}

Using the Bruhat decomposition one finds that
\[
\mu(\gln) = (\L^n - 1)(\L^n - \L)\ldots (\L^n - \L^{n-1})
\]
and hence that the  motive of $\gln$ is invertible in $\kvarhat$.
This will be important below.  For other groups we are interested in
we have:

\begin{proposition}\label{groupprop}
Let $G$ be a connected split semisimple group over $k$.  Then 
$$\mu(G)=\L^{\dim G}\prod_{i=1}^r(1-\L^{-d_i})$$
in $\kvarhat$. Here $r$ is the rank of $G$ and the $d_i$ are the
numbers one higher than the exponents of $G$. 
\end{proposition}
\begin{proof}
We choose a
Borel subgroup $B$ of $G$ with maximal torus $T$ and unipotent radical
$U$. Since $T$-bundles and $U$-bundles over varieties are
Zariski-locally trivial, we have $\mu(G)=\mu(G/P)\mu(T)\mu(U)$. The
torus $T$ is a product of multiplicative groups, so
$\mu(T)=(\L-1)^r$. The unipotent group $U$ is an iterated extension of
additive groups, so $\mu(U)=\L^u$, where $u={\frac{1}{2}(\dim G-r)}$
is the dimension of $U$. Finally, the
flag variety $G/B$ has a cell decomposition coming from the Bruhat
decomposition, and we have $\mu(G/B)=\sum_{w\in W}\L^{\ell(w)}$, where
$W$ is the Weyl group and $\ell(w)$ the length of a Weyl group
element. We have 
$$(\L-1)^r\sum_{w\in W}\L^{\ell(w)}=\prod_{i=1}^r(\L^{d_i}-1)\,,$$
by Page~150 of \cite{kane:00} or Page~155 of \cite{carter:72}, and
hence
$$\mu(G)
=\L^u\prod_{i=1}^r(\L^{d_i}-1)
=\L^{\dim G}\prod_{i=1}^r(1-\L^{-d_i})\,,$$
since $u+\sum_{i=1}^r d_i=\dim G$, by Solomon's theorem, see Page~320
of~\cite{kane:00}. 
\end{proof}

\subsection{The motive of an algebraic stack}

All our algebraic stacks will be Artin stacks, locally of finite
type, all of whose geometric stabilizers are linear algebraic
groups.  We will simply refer to such algebraic stacks as {\em
stacks with linear stabilizers}.

By a result of Kresch (Proposition~3.5.9 in \cite{kresch:99}), every
stack with linear stabilizers admits 
a stratification by locally closed substacks all of which are
quotients of a variety by $\gln$, for various $n$.  Note that
unless the stack $\sX$ is of finite type, there is no reason why
such a stratification should be finite.

Let us remark that for any stack $\sZ$, the reduced substack
$\sZ^\mathrm{red}\subset \sZ$ is locally closed, so that $\sX$ and
$\sX^\mathrm{red}$ have the same stratifications by locally closed
reduced substacks.

\begin{defn}
We call a  stack $\sX$ with linear stabilizers {\bf essentially
}of finite type, if it admits a countable stratification
$\sX=\bigcup \sZ_i$, where each $\sZ_i$ is of finite type and
$\dim \sZ_i\rightarrow -\infty$ as $i\to\infty$.
\end{defn}

Every stack with linear stabilizers which is essentially of finite
type admits countable stratifications $\sX=\bigcup \sZ_i$, where 
$$\lim_{i\to\infty}\dim\sZ_i=-\infty$$
and every $\sZ_i$ is a global quotient of a $k$-variety $X_i$ by a suitable
$\gln[n_i]$.   We call such stratifications {\em standard}. 

Let 
$$\sX=\bigcup_{i=0}^\infty[X_i/\gln[n_i]]$$
be  a standard stratification of the
essentially of finite type stack $\sX$. Define
$$\mu(\sX)=\sum_{i=0}^\infty\frac{\mu(X_i)}{\mu(\gln[n_i])}\,.$$
Note that the infinite sum converges in $\kvarhat$, by our assumptions.

The next lemma implies that our definition of $\mu(\sX)$, the motive
of the stack $\sX$, does not depend on the choice of a standard
stratification of $\sX$.

\begin{lemma}
Let $\sX\cong[X/\gln]$ be a global quotient stack, where $X$ is a
variety.  Let $\sX=\bigcup_{i=1}^{N}\sZ_i$ be a stratification of $\sX$
by locally closed substacks $\sZ_i$, which are, in turn, global
quotient stacks $\sZ_i\cong[X_i/\gln[n_i]]$. Then
$$\frac{\mu(X)}{\mu(\gln)}=\sum_{i=1}^N\frac{\mu(X_i)}{\mu(\gln[n_i])}$$
in $\kvarhat$. 
\end{lemma}
\begin{proof}
Let $Z_i$ be the preimage of $\sZ_i\subset \sX$ in $X$ under the
structure morphism $X\to\sX$.  Then $\sZ_i\cong[Z_i/\gln]$ and
$\sZ_i\cong[X_i/\gln[n_i]]$. Define $Y_i$ as the fibered product
$$\xymatrix{
Y_i\rto\dto & Z_i\dto\\
X_i\rto & \sZ_i}$$
Then $Y_i\to X_i$ is a principal $\gln$-bundle and $Y_i\to Z_i$ is a
principal $\gln[n_i]$-bundle.  Since $\gln$-bundles are always
Zariski-locally trivial, we conclude that $\mu(Y_i)=\mu(X_i)\mu(\gln)$
and $\mu(Y_i)=\mu(Z_i)\mu(\gln[n_i])$. Thus we have
$$
\frac{\mu(X)}{\mu(\gln)}= \sum_{i=1}^N\frac{\mu(Z_i)}{\mu(\gln)}
=\sum_{i=1}^N\frac{\mu(X_i)}{\mu(\gln[n_i])}
$$
as required.
\end{proof}

According \cite[Proposition 3.5.5]{kresch:99} the class
of algebraic stacks $\sX$ for which $\mu(\sX)$ makes sense
includes all Deligne-Mumford stacks of finite type.

\subsection{The torsor relations}

An essential ingredient in the definition of the motive of a stack
with linear stabilizers was the fact that every $\gln$-principal bundle
over a variety is Zariski locally trivial.  This implies that if $P\to X$ is a
principal $\gln$-bundle, then
\begin{equation}\label{torsor}
\mu(P)=\mu(X)\mu(\gln)\,,
\end{equation}
even if $X$ is a stack (where $\gln$-bundles are not necessarily
Zariski locally trivial any longer). 

In Section~\ref{p1}, we will need (\ref{torsor}) to hold for more
general groups than $\gln$. This is why we make the following
definition.

\begin{defn}
Fix an algebraic group $G$. 
We define $\kvarhatG$ to be the quotient of the ring $\kvarhat$ by
the ideal generated by all elements
$$\mu(P)-\mu(X)\mu(G)$$
where $X$ is a $k$-variety, and $P\to X$ is a $G$-torsor. 
\end{defn}

\begin{lemma}\label{torcon}
Let $\sX$ be an essentially of finite type stack with linear
stabilizers 
and $P\to\sX$ a $G$-torsor.  Then $P$ is also essentially of finite
type with linear stabilizers and we have
$$\mu(P)=\mu(\sX)\mu(G)$$
in $\kvarhatG$. 
\end{lemma}

\begin{example}\label{groupex}
Let $G$ be a connected split semisimple group over $k$. Then we have
$$\mu(BG)=\L^{-\dim G}\prod_{i=1}^r(1-\L^{-d_i})^{-1}\,$$
in $\kvarhatG$. Indeed,the torsor relation for $G$ (or Lemma
\ref{torcon}) implies that we 
have $\mu(BG)=\mu(G)^{-1}$. Now apply Proposition~\ref{groupprop}.
\end{example}

\begin{rmk}
Introducing the torsor relation $\mu(P)=\mu(X)\mu(G)$ for disconnected
$G$ kills $\kvarhat$.  For example, consider the $\mu_2$-torsor
$\G_m\to\G_m$. If $\Char k\not=2$, $\mu_2\cong \zz/2$ and the torsor
relation would imply $\L-1=2(\L-1)$ and hence $1=2$, as $\L-1$ is
invertible. 
\end{rmk}

\begin{rmk}
For connected $G$, the ring $\kvarhatG$ is non-trivial.  For
example, the $\ell$-adic 
Hodge-Poincar\'e characteristic (called the Serre characteristic by
some authors), factors through $\kvarhatG$. This follows from the
fact that a connected group cannot act non-trivially on its own
$\ell$-adic cohomology.
By the same token, the singular Hodge-Poincar\'e characteristic (in
case $k=\cc$) and the counting measure (in case $k=\ffq$) also factor
through $\kvarhatG$. 
\end{rmk}

\begin{rmk}
The second named author of this paper proves in the appendix that the
torsor relation (for split and connected linear algebraic
groups) holds in Voevodsky's category of effective geometrical
motives.  
\end{rmk}

\begin{rmk}
Recall that an algebraic  group $G$ over $k$ is called {\em special},
if all its torsors over $k$-varieties are Zariski-locally trivial. 
For special groups $G$, we have $\kvarhatG=\kvarhat$. Special groups
include $\sln$ and the symplectic groups $\mathop{\rm
  Sp}_{2n}$.

One may ask to what extent 
$\kvarhatG$ differs from  $\kvarhat$, for various groups $G$.
\end{rmk}

\section{The Main Conjecture}

Let $G$ be a split connected semisimple algebraic group over
$k$.
We denote by $d_1,d_2,\ldots,d_r$, where $r$ is the rank of $G$,  the
numbers one higher 
than the exponents  $G$. It will be important
below that $d_i\ge 2$.
Let $W$ be the
Weyl group of $G$ and $X(T)$ the character group of
a maximal torus $T$ of $G$. Then $W$ acts on the symmetric algebra
of $X(T)$. The $d_i$ are characterized by the fact that the ring of
invariants has generators in degrees 
$d_i$, see \cite{chevalley:55}.

Let $C$ be a smooth projective geometrically connected algebraic curve
over $k$, of genus $g$. Denote by $\sym{C}$ the $n$th symmetric power
of $C$. Recall that the motivic zeta function of $C$ is the power series
\[
\Zeta(C,u) = \sum_{n=0}^\infty \mu(\sym{C})u^n\quad\in\quad\kvarhat[[u]]\,.
\]
It is known that this function is in fact rational in $u$, see
\cite{kapranov:00} and \cite[\S 3]{larsen:02}. The denominator
is
\[
(1-u)(1-\L u)
\]
and hence evaluating the zeta function at $u=\L^{-n}$ makes sense
when $n\ge 2$.

We denote by $\bun{}{G,C}$ the moduli stack of $G$-torsors
over $C$.  The motive of $\bun{}{G,C}$ is defined by the following
lemma. 

\begin{lemma}
The stack $\bun{}{G,C}$ is essentially of finite type with linear
stabilizers.
\end{lemma}
\begin{proof} See \cite{behrend:t} or \cite{behrend:05} for full
details. For foundational results on the canonical parabolic
the reader is referred to \cite{behrend:95}.
The automorphism group scheme of a $G$-bundle $E$ is equal to the
scheme of global sections $\Gamma(C,\Aut G)$, where $\Aut(G)$ is the
group scheme over $X$ of automorphisms of $E$. Since $\Aut(G)$ an
affine over $C$ and $C$ is projective, $\Gamma(C,\Aut G)$ is affine,
hence linear.  Thus $\bun{}{G,C}$ has linear stabilizers.

Choose a Borel subgroup  $B$ of $G$ and call parabolic subgroups of $G$
containing $B$ {\em standard}.  Then every $G$-torsor $E$ over $C$ has a
canonical reduction of structure group $F$ to a uniquely determined
standard parabolic $P\subset G$. Thus $E=F\times_PG$.  The degree of
(the Lie algebra of) the group scheme $\Aut(F)={}^FP=F\times_{P,ad}P$
is called the {\em degree of instability }of $E$. It is a non-negative
integer (and 0 if and only if $E$ is semi-stable). Note that we allow
$G$ itself to be a parabolic subgroup in this context.

For every $m\geq0$, the substack $\bun{\leq m}{}\subset\bun{}{G,C}$ of
torsors of degree of instability less than or equal to $m$ is open in
$\bun{}{G,C}$ and of finite type.  The substack $\bun{m}{}$ of torsors
of degree of instability equal to $m$ is locally closed in
$\bun{}{G,C}$ and of dimension $\dim P (g-1)-m$, which is certainly
less than or equal to $\dim G (g-1)-m$, so tends to $-\infty$, as $m$
goes to $\infty$.
\end{proof}

We now come to our main conjecture. 

\begin{conjecture}\label{special}
If $G$ is simply connected, we have
$$\mu(\bun{}{G,C})= \L^{(g-1)\dim G}
\prod_{i=1}^r \Zeta(C,\L^{-d_i})
$$
in  $\kvarhatG$. 
\end{conjecture}

\begin{remark}
The conjecture makes sense inside the ring $\kvarhat$, but we dare not
conjecture its truth in the absence of the torsor
relations for $G$.  The proof in the case of $C=\pp^1$ uses the torsor
relations in an essential way, as we use the formula
$\mu(BP)\mu(G)=\mu(G/P)$, for all parabolic subgroups $P$ of $G$. But
note that this requires the torsor relations only for the group $G$ and
no others.

Note also, how the formula in Example~\ref{groupex} can be thought of
as an analogue of our conjecture for $C$ replaced by $\spec
k$. Example~\ref{groupex} also relies on the torsor relation for $G$.
\end{remark}

We can generalize the conjecture to arbitrary split connected
semisimple $G$:

\begin{conjecture}\label{general}
We have 
$$\mu(\bun{}{G,C})=|\pi_1(G)|\,\L^{(g-1)\dim G}
\prod_{i=1}^r \Zeta(C,\L^{-d_i})
$$
in $\kvarhatG$.
\end{conjecture}

Heuristically, the general case follows from the simply connected case
because we expect
$\bun{}{G,C}$ to have $|\pi_1(G)|$ connected components, all with
motive equal to the motive of $\bun{}{\widetilde{G},C}$, where
$\widetilde{G}$ is the universal covering group of $G$. 

The rest of this paper is devoted to providing evidence for our
conjecture.

\section{Evidence from Gauge Field Theory}\label{gauge}

In this section, $k=\cc$.
Denote by
\[
\chi_c:\cvarhat\rightarrow \Z((t^{-1}))
\]
the Poincar\'e characteristic. We will check that
Conjecture~\ref{special} holds after applying $\chi_c$ to both sides.

For a smooth $\cc$-variety
$X$ of dimension $n$, we have
\begin{align*}
\chi_c(\mu(X))&= \sum_{i,j}(-1)^j\dim W^iH^j_c(X,\cc)t^i\\
&=t^{2n}\sum_{i,j}(-1)^j\dim W^iH^j(X,\cc)t^{-i}\\
&=t^{2n}P_w(X,t^{-1})\,,
\end{align*}
by Poincar\'e duality, where $P_w(X,t)$ is the Poincar\'e polynomial
of $X$ using {\em weights}.

The cohomology of a finite type $\cc$-stack is endowed with a mixed Hodge
structure.  It is constructed via simplicial resolutions of the
stack.  Because every $\cc$-stack $\sX$, which is essentially of finite type
with linear stabilizers, can be exhausted by finite type open
substacks, the cohomology $H^n(\sX,\cc)$ of $\sX$ also carries a mixed
Hodge structure (for every $n$, the space $H^n(\sX,\cc)$ is equal to
the $n$-th cohomology of a sufficiently large finite type open
substack of $\sX$). Thus, $\sX$ has a Poincar\'e series
$$P_w(\sX,t)=\sum_{i,j}(-1)^j\dim W^i H^j(\sX,\cc)t^i\,.$$

\begin{lemma}
For every essentially of finite type $\cc$-stack with linear stabilizers
$\sX$ which is smooth, we have
$$\chi_c(\mu(\sX))=t^{2\dim\sX}P_w(\sX,t^{-1})\,.$$
\end{lemma}
\begin{proof}
If $\sX=[X/\gln]$ is a global quotient, the formula holds by the Leray
spectral sequence for the projection $X\to\sX$.  

Suppose $\sX$ is smooth of finite type and $\sZ$ a smooth closed
substack.
Then $P_w(\sX,t)=P_w(\sX-\sZ,t)+t^{2\codim(\sZ,\sX)}P_w(\sZ,t)$.  This
follows easily from the scheme case by using a simplicial resolution
$X_\bullet $ of $\sX$ and the fact that
$P_w(\sX,t)=\sum_j(-1)^jP_w(X_j,t)t^j$. 

Putting these two remarks together, we get the lemma in the finite
type case. For the general case, we choose a stratification
$\sX=\bigcup_{i=0}^\infty\sZ_i$, such that every
$\sX_n=\bigcup_{i=0}^n\sZ_i$ is a finite type open substack of $\sX$ and
$\lim_{i\to\infty}\dim\sZ_i=-\infty$.  Then we have
\begin{align*}
\chi_c(\mu(\sX))& =\chi_c(\lim_{n\to\infty}\mu(\sX_n))\\
&=\lim_{n\to\infty}\chi_c(\mu(\sX_n))\\
&=t^{2\dim\sX}\lim_{n\to\infty}P_w(\sX_n,t^{-1})\\
&=t^{2\dim\sX}P_w(\sX,t^{-1})\,,
\end{align*}  
where the last equality follows from the fact that for fixed $p$, the
cohomology group $H^p(\sX_n,\qq)$ stabilizes, as $n\to\infty$. 
\end{proof}

From \cite{macdonald:62} we have that
\[
\chi_c(\Zeta(C,u)) = \frac{(1+ut)^{2g}}{(1-u)(1-ut^2)}\,.
\]
So because $\bun{}{G,C}$ is smooth of dimension $\dim G(g-1)$,
our conjecture for the simply connected case becomes
\begin{equation}\label{simply}
P_w(\bun{}{G,C},t) = \prod_{i=1}^r
\frac{(1+t^{2d_i -1})^{2g}}{(1-t^{2d_i})(1-t^{2(d_i-1)})}.
\end{equation}
upon applying $\chi_c$ to both sides.

The Hodge structure on the cohomology of $\bun{}{G,C}$ has been
computed by Teleman~\cite{teleman:98}. In fact, Teleman shows
(Proposition~(4.4) of [ibid.]) that the
Hodge structure on 
$H^\ast(\bun{}{G,C})$ is pure, i.e., that the Poincar\'e series $P_w$
using weights is equal to the Poincar\'e series $P$ using Betti
numbers.  Thus we are reduced to computing Betti numbers of
$\bun{}{G,C}$. Atiyah and Bott \cite{atiyah:82} show that
$$
H^\ast(\bun{}{G,C})
\cong
H^\ast(G)^{\otimes 2g} \otimes H^\ast(BG)\otimes H^\ast(\Omega
G)\,.
$$

It is well known, see for example \cite{borel:53}, that we have the
following formulas for Poincar\'e series:
\begin{eqnarray*}
  P(G,t) &=& \prod_{i=1}^r (1 + t^{2d_i -1}) \\
  P(BG,t)&=& \prod_{i=1}^r \frac{1}{1-t^{2d_i}}\,.
\end{eqnarray*}
For the loop group $\Omega G$ we have, see \cite{bott:56} or \cite{garland:75},
\[
P(\Omega G,t) = \prod_{i=1}^r \frac{1}{1- t^{2(d_i -1)}}\,.
\]
This proves the desired formula~(\ref{simply}).

\begin{remark}
With no extra effort we can generalize the results of this section to
the  Hodge-Poincar\'e or Serre characteristic.
Recall that the Serre
characteristic $s(X;u,v)$ of a $\cc$-variety $X$ is defined
as
$$s(X;u,v)=\sum_{i,p,q}(-1)^i h^{p,q}H^i(X,\cc)u^pv^q$$
The Serre characteristic is also well-defined for elements of
$\cvarhat$ and for 
essentially of finite type $\cc$-stacks. If we apply the Serre
characteristic to our conjecture (in the simply connected case) we
obtain
\begin{equation}
s(\bun{}{G,C};u,v)=\prod_{i=1}^r
\frac{(1+u^{d_i}v^{d_i-1})^g(1+u^{d_i-1}v^{d_i})^g} 
{(1-u^{d_i}v^{d_i})(1-u^{d_i-1}v^{d_i-1})}
\end{equation}
This is exactly what Teleman proves in Proposition~(4.4) of~\cite{teleman:98}.
\end{remark}

\section{Evidence from Automorphic Forms}\label{auto}

In this section $k={\mathbb F}_q$.  The counting measure
$\#:\fqvar\to\zz$ extends to a ring morphism
$$\#:\fqvar[\L^{-1}]\longrightarrow \Q\,,$$
but this extension is not continuous, so there is no natural extension
of $\#$ to $\fqvarhat$ with values in $\R$. Still, we can make sense
of $\#$ on 
a certain  subring   of {\em convergent
}motives.

Choose a an embedding $\qqlbar\hookrightarrow\cc$. 
We have the compactly supported Frobenius characteristic 
$$F_c:\fqvarhat\longrightarrow \cc((t^{-1}))\,,$$
which is characterized by
$$F_c(\mu X,t)=\sum_{i,j}(-1)^j\tr F_q|W^iH^j_c(\ol{X},\qql)t^i\,,$$
for varieties $X$ over $\ffq$. Here $H^j_c(\ol{X},\qql)$ is the
$\ell$-adic \'etale cohomology with compact supports of the lift
$\ol{X}$ of $X$ to the algebraic closure of $\ffq$.  The (geometric)
Frobenius acting on $\ell$-adic cohomology is denoted by $F_q$. 

\begin{defn}
We call an element $x\in\fqvarhat$ with compactly supported Frobenius
characteristic $F_c(x,t)=\sum_{n}a_nt^{-n}$ {\bf convergent }if the
series $\sum_{n}a_n$ converges absolutely in $\cc$.
If this is the case, we call the sum
$\sum_n{a_n}$ the {\bf counting measure }of $x$, notation $\#(x)$.
\end{defn}

The convergent elements form a subring $\fqvarhat_{\rm conv}$ of
$\fqvarhat$, and we have a well-defined counting measure
$$\#:\fqvarhat_{\rm conv}\longrightarrow \cc\,,$$
which is a ring morphism. 
Note that $\#$ is not continuous.  For example, the sequence
$q^n/\L^n$ converges to zero in $\fqvarhat$, but  its counting
measure converges to~1.

\begin{lemma}
Every  finite type $\ffq$-stack with linear stabilizers $\sX$ has
convergent motive $\mu(\sX)$. Moreover, $\#(\mu \sX)$ is equal to
$\#\sX(\ffq)$, the 
number of rational points of $\sX$ over $\ffq$, counted in the stacky
sense, i.e., we count isomorphism classes of the category $\sX(\ffq)$,
weighted by the reciprocal of the number of automorphisms. 
\end{lemma}
\begin{proof}
This lemma reduces to the Lefschetz trace
formula for $F_q$ on the compactly supported cohomology of an
$\ffq$-variety.  The reduction uses the simple fact that
$\#[X/\gln](\ffq)=\#X(\ffq)/\#\gln(\ffq)$.
\end{proof}

Because of the non-continuity of the counting measure, this lemma does
not generalize to all essentially finite type stacks over $\ffq$.  But
we do have a result for certain smooth stacks:

\begin{lemma}\label{complicated}
Let $\sX$ be a smooth stack with linear stabilizers over $\ffq$.
Suppose that $\sX$ has a stratification
$\sX = \cup_{i=0}^\infty \sZ_i$ by smooth substacks $\sZ_i$, such that
for 
every $n$  the stack 
$\sX_n = \cup_{i=0}^n\sZ_i$ is an open substack of finite type 
and 
$$\sum_{n=0}^\infty q^{-\codim(\sZ_n,\sX)}\sum_{i,j} \dim W^i
H^j(\ol{\sZ}_n,\qql)q^{-i/2}<\infty\,.$$
Then $\sX$ is essentially of finite type, its
motive $\mu\sX$ is convergent and $\#(\mu\sX)=\#\sX(\ffq)$.
\end{lemma}
\begin{proof}
Let us emphasize that we assume that for
every $i$, the substack $\sZ_i$ is non-empty and its codimension
inside $\sX$ is constant. Let us  denote this codimension by $c_i$. 

Let us also remark that our assumptions imply that
$\sX$ is essentially of finite type and that 
$\lim_{n\to\infty} c_n=\infty$.  We may also assume, without loss of
generality, that the dimension of $\sX$ is constant. 

First, we will prove that the trace of the arithmetic Frobenius on the
$\ell$-adic cohomology of $\sX$ converges absolutely to
$q^{-\dim\sX}\#\sX(\ffq)$.

There is a spectral sequence of finite dimensional $\qql$-vector
spaces
$$E_1^{pq}=H^{p+q-2c_p}(\ol{\sZ}_p,\qql(-c_p))\Longrightarrow
H^{p+q}(\ol{\sX},\qql) \,.$$
Even though this is not a first quadrant spectral sequence, we do have
that for every $n$ there are only finitely many $(p,q)$ with $p+q=n$ and
$E_1^{pq}\not=0$, so this spectral sequence does converge. 

Our assumption on $\sX$ implies that the arithmetic Frobenius $\Phi_q$
acting on $E_1$ has absolutely convergent trace. Thus we get the same
result for this trace, no matter in which order we perform the
summation. Thus, using the trace formula for the arithmetic Frobenius
on finite type smooth stacks with linear stabilizers (see
\cite{behrend:93})  we have
\begin{align*}
\#\sX(\ffq)&= \sum_{p=0}^\infty\#\sZ_p(\ffq)\\
&= \sum_{p=0}^\infty q^{\dim\sZ_p}\tr\Phi_q|H^\ast(\ol{\sZ}_p,\qql)\\
&= q^{\dim\sX}\sum_{p=0}^\infty
\tr\Phi_q|H^\ast(\ol{\sZ}_p,\qql(-c_p))\\
&= q^{\dim\sX} \tr\Phi_q|H^\ast(\ol{\sX},\qql)
\end{align*}
In particular, we see that $\#\sX(\ffq)$ is finite.

Next we will examine the motive of $\sX$. Note that for smooth stacks
of finite type $\sY$, we have 
$$F_c(\mu\sY,t)=(qt^2)^{\dim\sY} \Phi(\sY,t^{-1})\,$$
where $\Phi$ is the Frobenius characteristic defined using the
arithmetic Frobenius acting on cohomology
without compact supports:
$$\Phi(\sY,t)=\sum_{i,j}(-1)^j\tr\Phi_q|W^iH^j(\ol{\sY},\qql)t^i\,.$$
This is essentially Poincar\'e duality for smooth varieties. 
Thus we have
\begin{align*}
F_c(\mu\sX,t)&=\lim_{n\to\infty}F_c(\mu\sX_n,t)\\
&=(qt^2)^{\dim\sX}\lim_{n\to\infty}\Phi(\sX_n,t^{-1})\\
&=(qt^2)^{\dim\sX}\Phi(\sX,t^{-1})\,.
\end{align*}
So to prove that $\mu\sX$ is convergent, we need to prove that 
$$\sum_i\Bigl\lvert\sum_j(-1)^j\tr\Phi_q|
W^iH^j(\ol{\sX},\qql)\Bigr\rvert<\infty\,.$$  
But our spectral sequence implies that
$$\sum\dim W^iH^j(\ol{\sX},\qql)q^{-i/2}<\infty\,,$$
which is a stronger statement. 
So we see that $\mu\sX$ is, indeed, convergent and its counting
measure takes the value
$$\#(\mu\sX)=q^{\dim\sX}\tr\Phi_q|H^\ast(\ol{\sX},\qql)\,.$$
This we have seen above to be equal to $\#\sX(\ffq)$. 
\end{proof}

We say that a morphism of stacks $\sZ\rightarrow \tilde{\sZ}$
is a universal homeomorphism if it is representable, finite,
surjective and radical.

\begin{lemma}\label{sufficient}
Lemma~\ref{complicated} is still valid if we only assume the morphisms
$\sZ_i\to \sX$ to be universal homeomorphisms onto their image. 
\end{lemma}
\begin{proof}
Let  $Z\to X$ be a morphism of finite type smooth
schemes which factors as $ Z\to\widetilde Z\to X$, where
$\pi:Z\to\widetilde Z$ is a universal homeomorphism and
$i:\widetilde Z\to X$ a closed immersion with complement $U$. We
have a long exact sequence
$$\ldots\to H^\ast(\widetilde Z,i^!\qql)\to 
H^\ast( X,\qql)\to 
H^\ast( U,\qql)\to\ldots
$$
Let $c=\dim X-\dim Z$.  We have
$$H^{\ast-2c}(Z,\qql(-c))=H^\ast(Z,\pi^!i^!\qql)\,$$
because $Z$ and $X$ are smooth.  Now pulling back via $\pi$ induces an
isomorphism of \'etale sites (see \cite[Expose IX,4.10]{SGA1}).
As $\pi_\ast$ is  the
right adjoint of $\pi^\ast$,  it is the inverse of $\pi^\ast$ and hence also
a left adjoint of $\pi^\ast$.  Since $\pi$ is proper, 
we conclude that $\pi^!=\pi^\ast$. Thus, we have
$$H^\ast(Z,\pi^!i^!\qql)=H^\ast(Z,\pi^\ast i^!\qql)=H^\ast(\widetilde
Z,i^!\qql)\,,$$
Thus we have a natural long exact sequence
$$\ldots\to H^{\ast-2c}( Z,\qql(-c))\to 
H^\ast( X,\qql)\to 
H^\ast( U,\qql)\to\ldots
$$
This result extends to stacks and filtrations of schemes and stacks
consisting of more than two pieces.
\end{proof}

\begin{lemma}
The motive of $\bun{}{G,C}$ is convergent. Moreover,
$\#\mu(\bun{}{G,C})=\#\bun{}{G,C}(\ffq)$. 
\end{lemma}
\begin{proof}
The hypotheses of Lemma~\ref{complicated}, or rather 
its generalization~\ref{sufficient}, are satisfied by the stack
$\bun{}{G,C}$. We may consider the strata $\bun{P,m}{G,C}$, which
contain the bundles $E$ which canonically reduce to the standard
parabolic $P$ of $G$ and whose degree of instability is equal to
$m$, see \cite{behrend:95}. 
These strata are not known to be smooth, but the canonical
morphism $\bun{{ ss},m}{P,C}\to\bun{P,m}{G,C}$ is a universal
homeomorphism.    Here  
$\bun{{  ss},m}{P,C}$ is the open substack of 
$\bun{}{P,C}$ consisting of semi-stable bundles of positive
(multi-)degree, giving rise to degree of instability $m$ when
extending the group to $G$. 

If $H$ is the quotient of $P$ by its unipotent radical, the induced
morphism $\bun{ss,m}{P,C}\to\bun{ss,m}{H,C}$ induces an isomorphism on
$\ell$-adic cohomology, because it is an iterated torsor for  vector bundle
stacks. 

This leaves us with proving the convergence of
$$\sum_P\sum_{m=1}^{\infty}q^{-m+(1-g)\dim R_uP}\sum_{i,j}\dim W^i
H^j(\ol{\bun{}{}}_{H,C}^{ss,m},\qql)q^{-i/2}\,.$$
This is not difficult to do using  the fact that for fixed $H$, all
$\bun{ss,m}{H,C}$ are isomorphic to a finite set among them. 
\end{proof}

By this lemma, both sides of our conjectured formula are in the
subring $\fqvarhat_{\rm conv}$.  We can thus apply the counting
measure $\#$ to our conjecture.  Doing this we obtain:
\begin{equation}\label{zeta}
\#\bun{}{G,C}(\ffq)=|\pi_1(G)|q^{(g-1)\dim G}\prod_{i=1}^r \zeta_K(d_i)
\end{equation}
Here $\zeta_K(s)$ is the usual zeta function of the function field $K$ of
the curve $C$ over $\ffq$.  It is obtained from the motivic zeta
function $Z(C,u)$ of the curve $C$ by applying the counting measure
and making the substitution $u=q^{-s}$. 

Formula~(\ref{zeta}) is classical, at least in the simply connected
case.  Let us recall how it is proved. 
We consider the ad\`ele ring $\A_K$ of the global field $K$ and 
notice that the groupoid $\bun{}{G,C}(\ffq)$ is equivalent to the
transformation groupoid of the action of $G(K)$ on $G(\A_K)/\kk$,
where $\kk=\prod_{P\in C}G(\widehat{\O}_{C,P})$ is the canonical
maximal compact subgroup of $G(\A_K)$. 

The transformation groupoid of the $G(K)$-action on $G(\A_K)/\kk$ is
equivalent to the transformation groupoid of the $\kk$-action on
$G(K)\backslash G(\A_K)$. The groupoid number of points of the latter
transformation groupoid can be calculated as 
\begin{equation}\label{tau}
\frac{\vol(G(K)\backslash G(\A_K))}{\vol{\kk}}\,,
\end{equation}
where $\vol$ denotes any Haar measure on the locally compact group
$G(\A_K)$.  This is  a simple measure theoretic argument using 
$\sigma$-additivity.  

There is a standard normalization of the Haar measure on $G(\A_K)$
known as the Tamagawa measure.  With respect to this measure the
numerator of (\ref{tau}) is known as the Tamagawa number of $G$,
notation $\tau(G)$. We conclude that
$$\#\bun{}{G,C}(\ffq)=\tau(G)\vol(\kk)^{-1}\,.$$
The volume of the maximal compact $\kk$ with respect to the Tamagawa
measure is easily calculated.  We get
$$\vol(\kk)=q^{(1-g)\dim G}\prod_{i=1}^r\zeta_K(d_i)^{-1}\,,$$
see \cite{behrend:05}, and thus
$$
\#\bun{}{G,C}(\ffq)=\tau(G)q^{(g-1)\dim G}\prod_{i=1}^r \zeta_K(d_i)\,.
$$
Comparing this with our conjecture (\ref{zeta}) we see that the
conjecture becomes equivalent to 
\begin{equation}\label{tamagawa}
\tau(G)=|\pi_1(G)|\,.
\end{equation}
In the simply connected case, the fact $\tau(G)=1$ was 
proved by Harder \cite{harder:74}. 
The results of \cite{ono:65}
remain true in the function field case (see~\cite{behrend:05})
and from these it follows that the Tamagawa number
of a general connected split semisimple group $G$ is equal to
$|\pi_1(G)|$.

\section{The Case of ${\rm Sl}_n$}\label{sln}

In the section we prove our conjecture in the case
where the group is $G=\sln$. 
Recall that the exponents of
$\sln$ are $2,3,\ldots , n$.  Thus our conjecture states that
$$\mu(\bun{}{\sln,C})=\L^{(n^2-1)(g-1)}\prod_{i=2}^n Z(C,\L^{-i})\,.$$

To calculate the motive of
$\bun{}{\sln}$, note that the inclusion $\sln\hookrightarrow\gln$
defines a morphism of stacks $\bun{}{\sln,C}\to\bun{}{\gln,C}$, whose
image is a smooth closed substack $\bun{}{\det}$ of
$\bun{}{\gln}$. Moreover, $\bun{}{\sln}$ is a $\Gm$-bundle over
$\bun{}{\det}$. Thus we have
$$\mu(\bun{}{\sln})=(\L-1)\mu(\bun{}{\det})\,.$$
We can interpret $\bun{}{\det}$ is the stack of vector bundles over
$C$ with trivial determinant. 

We will use the construction of matrix divisors in
\cite{bifet:94}. Let $D$ be an effective divisor on $C$.  
We denote by $\div(D)$ the Quot scheme parameterizing subsheaves
$$E\hookrightarrow \O_C(D)^n\,,$$
where $E$ is a locally free sheaf of rank $n$ and degree 0 on $C$. 
The scheme $\div(D)$ is smooth and proper of dimension $n^2\deg D$. 

Let $\div_{\det}(D)\subset\div(D)$ by the closed subscheme defined by
requiring the determinant of $E$ to be trivial.  This is a smooth
subscheme of codimension $g$.  (See~\cite{dhillon:05} for the
proof of this.) 

Now let us fix, for the moment, an integer $m\geq0$ and consider the
finite type open substack $\bun{\leq m}{\det}$, of bundles whose
degree of instability is at most $m$.  Let $D$ be an effective divisor
of sufficiently high degree, such that $H^1(E,\O(D)^n)=0$, for all
bundles $E$ in $\bun{\leq m}{\det}$. Then the vector spaces 
$\Hom(E,\O(D)^n)$, for $E\in \bun{\leq m}{\det}$, are the fibres of a
vector bundle $W^{\leq  m}(D)$ over $\bun{\leq m}{\det}$. The rank of
this vector bundle is $n^2(\deg D +1-g)$. 

Let $W_0^{\leq m}(D)\subset W^{\leq m}(D)$ be the open locus of
injective maps $E\to \O(D)^n$. Note that 
$$W_0^{\leq m}(D)=\div^{\leq m}_{\det}(D)$$
is the open subvariety of $\div_{\det}(D)$ parameterizing subsheaves
$E\subset \O(D)^n$ of degree of instability at most $m$. 
$$\xymatrix{
W^{\leq m}(D)\drto_{\text{vector bundle}}&
W_0^{\leq m}(D) \ar@{=}[r]\ar@{_{(}->}[l]\dto & \div^{\leq m}_{\det}(D)
 \ar@{^{(}->}[r]\dlto & \div_{\det}(D)\dlto\\
&\bun{\leq m}{\det}\ar@{^{(}->}[r] & \bun{}{\det}}$$

\begin{lemma}
Let $E$ and $F$ be vector bundles of equal rank on $C$.  Let $D$ be an
effective divisor on $C$ such that $H^1(E,F(D))$ vanishes.  Then the
locus of the non-injective maps inside $\Hom(E,F)$ has codimension at
least $\deg D$. 
\end{lemma}
\begin{proof}
This is proved in Lemma~8.2 of~\cite{bifet:94}.  
\end{proof}

This lemma implies that
$$\lim_{\deg D\to\infty}\frac{\mu( W^{\leq m}_0)}{\L^{n^2(\deg D+1-g)}}=
\lim_{\deg D\to\infty}\frac{\mu (W^{\leq m})}{\L^{n^2(\deg D+1-g)}}$$
inside $\kvarhat$. Thus we have
\begin{align*}
\mu(\bun{}{\det}) 
& = \lim_{m\to\infty}\mu(\bun{\leq m}{\det})\\
& = \lim_{m\to\infty}\,\lim_{\deg D\to\infty}\frac{\mu(W^{\leq
    m}(D))}{\L^{n^2(\deg D+1-g)}}\\ 
& = \lim_{m\to\infty}\,\lim_{\deg D\to\infty}\frac{\mu(\div^{\leq
    m}_{\det}(D))}{\L^{n^2(\deg D+1-g)}}\\ 
& = \lim_{\deg D\to\infty}\,\lim_{m\to\infty}\frac{\mu(\div^{\leq
    m}_{\det}(D))}{\L^{n^2(\deg D+1-g)}}\\ 
& = \lim_{\deg D\to\infty}\frac{\mu(\div_{\det}(D))}{\L^{n^2(\deg
    D+1-g)}}\,.
\end{align*}
Therefore, the conjecture translates into
$$\lim_{\deg D\to\infty}\frac{\mu(\div_{\det}(D))}{\L^{n^2(\deg
    D+1-g)}}=\frac{\L^{(n^2-1)(g-1)}}{\L-1}\prod_{i=2}^n Z(C,\L^{-i})$$
or, in other words,
\begin{multline}\label{limeq}
\lim_{\deg D\to\infty}\frac{\mu(\div_{\det}(D))}{\L^{n^2(\deg
    D+1-g)}}\\=\frac{\L^{(n^2-1)(g-1)}}{\L-1}
\sum_{\mathbf{m}=(m_2,\ldots,m_n)}\mu(C^{(\mathbf{m})})\L^{-\sum_{i=2}^n
    i m_i}\,.
\end{multline}
Here the sum ranges over all $(n-1)$-tuples of non-negative integers
and we use the abbreviation
$C^{(\mathbf{m})}=C^{(m_2)}\times\ldots\times C^{(m_n)}$.

It remains to calculate the motive of $\div_{\det}(D)$. This we will do
by using the stratification 
induced by a suitable $\Gm$-action via the results of Bia\l
ynick-Birula~\cite{bb:73}.  Note that 
we can neglect strata whose codimension goes to infinity, as
$\deg D$ goes to infinity.  

Consider the action of the torus $\Gm^n$ on $\div(D)$ induced by the
canonical action on the vector bundle $\O_C(D)^n$. It restricts to an
action of $\Gm^n$ on $\div_{\det}(D)$. 

The fixed points of $\Gm^n$ on $\div(D)$ correspond to inclusions of
the form
$$\textstyle\bigoplus_{i=1}^n \O_C(D-E_i)\hookrightarrow\O_C(D)^n\,,$$
where $E_1,\ldots,E_n$ are effective divisors with $\sum \deg
E_i=n\deg D$ (see \cite{bifet:94}). Thus, the components of the fixed
locus in $\div(D)$ are indexed by ordered partitions
$\mathbf{m}'=(m_1,\ldots,m_n)$ 
of $n\deg D$ and the
component indexed by $\mathbf{m}'$ is isomorphic to 
$$C^{(\mathbf{m}')}=C^{(m_1)}\times\ldots\times C^{(m_n)}\,.$$

The intersection of the fixed component $C^{(\mathbf{m}')}$ with the
subvariety $\div_{\det}(D)$ is given by the condition that $\sum E_i$
be linearly equivalent to $nD$. Thus, if $m_1> 2g-2$, 
this intersection is a projective space bundle with fibre
$\pp^{m_1-g}$ over $C^{(\mathbf{m})}$, where
$\mathbf{m}=(m_2,\ldots,m_n)$.  So the motive of the fixed component
of $\div_{\det}(D)$ indexed by $\mathbf{m}'$ is given by 
$$\frac{\L^{m_1-g+1}-1}{\L-1}\mu(C^{(\mathbf{m})})\,.$$
We will see below, that we can neglect
the fixed components indexed by $\mathbf{m}'$ with $m_1\leq 2g-2$.

Now, consider the $\Gm$-action induced by the one-parameter subgroup
$\Gm\to \Gm^n$ given by
$$t\longmapsto(t^{\lambda_1},\ldots,t^{\lambda_n})\,,$$
where $(\lambda_1,\ldots,\lambda_n)$ is any strictly increasing
sequence of integers $\lambda_1<\ldots<\lambda_n$. The fixed locus of
$\Gm$ on $\div(D)$ is then the same as that of the whole torus
$\Gm^n$. We will study the strata
$$X^+_{\mathbf{m}'}=\{x\in\div(D)\mid\lim_{t\to0}tx\in
C^{(\mathbf{m}')}\}\,.$$  
and 
$$Y^+_{\mathbf{m}'}=\{x\in\div_{\det}(D)\mid\lim_{t\to0}tx\in
C^{(\mathbf{m}')}\cap\div_{\det}(D)\}\,.$$  
There is a morphism $X^+_{\mathbf{m}'}\to C^{(\mathbf{m}')}$ making
$X^+_{\mathbf{m}'}$ into a Zariski locally trivial affine  space
bundle over $C^{(\mathbf{m}')}$, see \cite{bb:73}. The rank of this
fibration is the same as the rank of the subbundle $N^+$ of $N$ on
which $\Gm$ acts with 
positive weights. Here $N$ is the normal bundle of
$C^{(\mathbf{m}')}$ inside $\div(D)$. 

The tangent space inside $\div(D)$ to the fixed point $P$ given by
$(E_1,\ldots,E_n)\in C^{(\mathbf{m}')}$ is equal to
$\bigoplus_{i,j}\Hom(\O_C(D-E_i), \O_{E_j})$ and the torus $\Gm^n$
acts on the summand $\Hom(\O_C(D-E_i),\O_{E_j})$ through the character
$\chi_i-\chi_j$, where $\chi_i$ is given by the $i$-th
projection $\chi_i:\Gm^n\to \Gm$. Thus we see that the fibre of $N^+$
over $P$ is
$$N^+_P=\textstyle\bigoplus_{i>j}\Hom(\O_C(D-E_i),\O_{E_j})\,,$$
and so the rank of $N^+$ is equal to $\sum_{i=1}^n(n-i)m_i$. 

If $P$ is in the subvariety $\div_{\det}(D)$, the tangent space to $P$
inside $\div_{\det}(D)$ is the kernel of the diagonal part of the
boundary map
$$\textstyle\bigoplus_{i,j}\Hom(\O_C(D-E_i), \O_{E_j})\longrightarrow
\bigoplus_{i,j}H^1(C,\O_C(E_i-E_j))$$
coming from the universal exact sequence
$$\textstyle0\longrightarrow \bigoplus_{i=1}^n\O_C(D-E_i)\longrightarrow
\O_C(D)^n\longrightarrow\bigoplus_{i=1}^n\O_{E_i}\longrightarrow 0\,.$$
(This is proved in \cite{dhillon:05}.) It follows that the rank of the
fibration 
$$Y^+_{\mathbf{m}'}\longrightarrow C^{(\mathbf{m}')}\cap
\div_{\det}(D)$$
is equal to $\sum_{i=1}^n(n-i)m_i$ as well. 

Now we see that the biggest stratum corresponds to an index
$\mathbf{m}'$ where $m_1$ attains the maximal value $n\deg D$.  The
dimension of all strata coming from $X_{\mathbf{m}'}$ or
$Y_{\mathbf{m}'}$ with $m_1\leq 2g-2$ is therefore bounded from above
by 
\begin{multline*}
\dim C^{(\mathbf{m}')}+(n-1)(2g-2)+(n-2)(n\deg D-(2g-2))\\=
n(n-1)\deg D+2g-2\,.
\end{multline*}
Hence their codimension inside $\div_{\det}(D)$ is bounded from below
by 
$$n^2\deg D-g -n(n-1)\deg D -(2g-2)=n\deg D-3g+2\,$$
which, indeed, goes to infinity with $\deg D$. We conclude that, up to
terms we are going to neglect, we have
$$\mu(\div_{\det}(D))\approx\sum_{\mathbf{m}'}
\frac{\L^{m_1-g+1}-1}{\L-1}\mu(C^{(\mathbf{m})})\L^{\sum_{i=1}^m(n-i)m_i}\,,$$
the sum ranging over all $\mathbf{m}'=(m_1,\ldots,m_n)$ with
$\sum_{i=1}^nm_i=n\deg D$.  We can rewrite this as
\begin{multline*}
\frac{\mu(\div_{\det}(D))}{\L^{n^2(\deg D+1-g)}}\\ \approx
\sum_{\mathbf{m}}\frac{\L^{-\sum_{i=2}^nm_i}-\L^{-n\deg D+g-1}}{\L-1}
\mu(C^{(\mathbf{m})})\L^{(n^2-1)(g-1)+\sum_{i=2}^n(1-i)m_i}\,
\end{multline*}
where the sum ranges over all $\mathbf{m}=(m_2,\ldots,m_n)$ with
$\sum_{i=2}^nm_i\leq n\deg D-2g+2$.  

As $\deg D$ goes to infinity this becomes an equality, in fact,
Equation~(\ref{limeq}), which  we set out to prove.

\section{The Case of ${\mathbb P}^1$}\label{p1}

In this section we use the Grothendieck-Harder classification of
torsors on ${\mathbb P}^1$ to prove the conjecture in the special
case that $C=\pp^1$.

We fix a split maximal torus $T$ inside $G$ and let $W$
be the Weyl group.   Let $\rx^*(T)$ (resp.\ $\rx_*(T)$) be
the character (resp.\ cocharacter) lattice.  We have the root system
$\Phi\subset \rx^\ast(T)$ and its dual $\Phi^\vee\subset\rx_\ast(T)$.

We also choose a Borel subgroup $B$
containing $T$. It determines bases $\Delta$ of $\Phi$  and
$\Delta^\vee$ of $\Phi^\vee$. 
Denote by $\rx_\ast(T)_\dom$ the dominant cocharacters
with respect to  $B$.  Recall that $\lambda\in \rx_\ast(T)$ is
dominant if and only if $(\lambda,\alpha)\geq0$, for all
$\alpha\in\Delta$. 
The set $\rx_\ast(T)_\dom$ is partially ordered:
$\lambda_1\leq\lambda_2$ if and only if $\lambda_2-\lambda_1$ is a
positive integral linear combination of elements of $\Delta^\vee$.

For a dominant cocharacter $\lambda\in\rx_\ast(T)_\dom$, denote by 
$$E_\lambda=\O(1)\times_{\Gm,\lambda} G$$
the $G$-bundle associated to the $\Gm$-bundle $\O(1)$ via the
homomorphism $\lambda:\Gm\to G$. (Think of the line bundle $\O(1)$ as
a $\Gm$-bundle.)

\begin{proposition}\label{sevenone}
Every $G$-bundle over $\pp^1_K$, for a field $K/k$, becomes isomorphic
to $E_\lambda$, for a unique $\lambda\in\rx_\ast(T)_\dom$,  after
lifting it to the algebraic closure of $K$. 
\end{proposition}
\begin{proof}
This result is obtained by combining the Grothendieck-Harder
classification of Zariski-locally trivial $G$-torsors by
$\rx_\ast(T)_\dom$ with the theorem of Steinberg, to the effect that
on $\pp^1$ over an algebraically closed field, all $G$-torsors are
Zariski-locally trivial.  See also Theorem~4.2 and Proposition~4.3
of~\cite{ramanathan:83}.
\end{proof}

By Proposition~\ref{sevenone}, the bundles $E_\lambda$, for
$\lambda\in\rx_\ast(T)_\dom$, give a complete set of representatives
for the points of the stack $\bun{}{G,\pp^1}$. Hence every point of
$\bun{}{G,\pp^1}$ is $k$-rational and its residual gerbe is trivial,
equal to $B\Aut E_\lambda$.

Recall that for $\sX$, a locally of finite type algebraic stack over
$k$ with set of points $|\sX|$, there is a topology on $|\sX|$, the
Zariski topology, such that open substacks of $\sX$ are in bijection
to open subsets of $|\sX|$.

Let us identify $|\bun{}{G,\pp^1}|$ with
$\rx_\ast(T)_\dom$. 

\begin{proposition}[Ramanathan]
Let $\lambda\in \rx_\ast(T)_\dom$.  Then the set of all
$\mu\in\rx_\ast(T)_\dom$ with $\mu\leq\lambda$ is open in the Zariski
topology on $|\bun{}{G,\pp^1}|$.
\end{proposition}
\begin{proof}
This is the content of Theorem~7.4 in~\cite{ramanathan:83}.
\end{proof}

It follows from this that the substack of $\bun{}{G,\pp^1}$ of torsors
isomorphic to $E_\lambda$ is locally closed.  Moreover, this substack
is necessarily equal to the substack $B\Aut E_\lambda$, because a
monomorphism of reduced algebraic stacks which is surjective on points
is an isomorphism.
Thus we have that
$$\bun{}{G,\pp^1}=\bigcup_{\lambda\in\rx_\ast(T)_\dom} B\Aut
E_\lambda$$ is a stratification of $\bun{}{G,\pp^1}$. 

To calculate the motive of $B\Aut E_\lambda$, fix the dominant
cocharacter $\lambda$.  Denote by $P$ the parabolic subgroup
of $G$ defined by $\lambda$ and by $U$ its unipotent
radical. The group $P$ is generated by $T$ and all root
groups $U_\alpha$, $\alpha\in\Phi$, such that
$(\lambda,\alpha)\geq0$. The group $U$ is generated by the $U_\alpha$
with $(\lambda,\alpha)>0$.  We will also use the Levi subgroup
$H\subset P$.  Note that $P=H\ltimes U$.

Via $\lambda:\Gm\to G$ the multiplicative group acts by conjugation on
$G$, $P$ and $U$.  We can use this action to twist $G$, $P$ and $U$ by
the $\Gm$-torsor $\O(1)$.  We denote the associated twisted groups by
$G_\lambda$, $P_\lambda$ and $U_\lambda$.  For example, 
$G_\lambda=\O(1)\times_{\Gm,\lambda,Ad}G$.

\begin{proposition}
We have $\Aut
E_\lambda=H\ltimes
\Gamma(\pp^1,U_\lambda)$. 
\end{proposition}
\begin{proof}
We have $\Aut E_\lambda=\Gamma(\pp^1,G_\lambda)=\Gamma(\pp^1,P_\lambda)=H\ltimes
\Gamma(\pp^1,U_\lambda)$. For more details, see~\cite{ramanathan:83},
Proposition~5.2.
\end{proof}

Note that for a semidirect product of linear algebraic groups $N$, $H$
we have $\mu B(H\ltimes N)=\mu(BH)\mu(BN)$. Thus, we may calculate
\begin{align*}
\mu (B\Aut E_\lambda) & = \mu B\big(
H\ltimes\Gamma(\pp^1,U_\lambda)\big)\\
&= \mu(BH) \mu \big(B\Gamma(\pp^1,U_\lambda)\big)\\
&=\frac{\mu (BP)}{\mu (BU)}{\mu \big(B\Gamma(\pp^1,U_\lambda)\big)}\\
&=\mu (BP)\frac{\mu (U)}{\mu\Gamma(\pp^1,U_\lambda)}\\
&=\frac{\mu(G/P)}{\mu G}\frac{\mu(U)}{\mu\Gamma(\pp^1,U_\lambda)}
\end{align*}
In the last equation we used the torsor relation for $G$, i.e.,
Lemma~\ref{torcon}.

Now let $\uU$ be the Lie algebra of $U$. We have
$\uU=\bigoplus_{(\lambda,\alpha)>0} \uU_\alpha$, where $\uU_\alpha$ is
the Lie algebra of $U_\alpha$. Since $T$ acts on each $\uU_\alpha$, we
obtain line bundles 
$$(\uU_\alpha)_\lambda=\O(1)\times_{\Gm,\lambda}\uU_\alpha\,.$$  
Note that the degree of $(\uU_\alpha)_\lambda$ is equal to
$(\lambda,\alpha)$. The 
group scheme $U_\lambda$ over $\pp^1$ is a successive extension of the
$(\uU_\alpha)_\lambda$, and therefore we have
$$\mu\Gamma(\pp^1,U_\lambda)
=\prod_{(\lambda,\alpha)>0}\L^{(\lambda,\alpha)+1}\,,$$
by Riemann-Roch and hence
$$\frac{\mu(U)}{\mu\Gamma(\pp^1,U_\lambda)}=\L^{-(\lambda,2\rho)}\,,$$
where $\rho$ is half the sum of all positive roots.

Denote by $W(\lambda)\subset W$ the subgroup of the Weyl group
generated by the reflections coming from simple roots orthogonal to
$\lambda$. The Bruhat decomposition for $G/P$ implies
$$\mu(G/P) = \sum_{w\in W/W(\lambda)} \L^{\ell(w)}\,,$$
where $\ell(w)$ is the minimum of all lengths in the coset
$wW(\lambda)$.

This finishes the analysis of the motive of $B\Aut E_\lambda$.
Putting everything together, we find:
\begin{align*}
\mu(\bun{}{G,\pp^1})&=
\sum_{\lambda\in\rx_\ast(T)_\dom} \mu(B\Aut E_\lambda)\\
&= \frac{1}{\mu
  G}\sum_{\lambda\in\rx_\ast(T)_\dom}\L^{-(\lambda,2\rho)} 
\sum_{w\in W/W(\lambda)} \L^{\ell(w)}\,.
\end{align*}

The combinatorics of summing the various powers of $\L$ is 
contained in \cite{kaiser:00}.  In fact, it is proved in [ibid.] that
$$\sum_{\lambda\in\rx_\ast(T)_\dom}\L^{-(\lambda,2\rho)} 
\sum_{w\in W/W(\lambda)} \L^{\ell(w)}=
|\pi_1(G)|\frac{P(W_{\rm aff},\L^{-1})}{P(W,\L^{-1})}\,.$$
Here the series $P(W_{\rm aff},t)$ (resp. $P(W,t)$)
is the Poincar\'e series of the affine Weyl (resp. Weyl) group. It is 
defined by
\[
P(W_{\rm aff},t) = \sum_{w\in W_{\rm aff}} t^{\ell(w)}.
\]
It is a result of Bott and Steinberg that
\[
\frac{P(W_{\rm aff},t)}{P(W,t)} = \prod_{i=1}^r (1-t^{d_i-1})^{-1}.
\]
Thus we may complete the calculation 
\begin{align*}
\mu(\bun{}{G,\pp^1})&=|\pi_1(G)|
\frac{1}{\mu
  G}\prod_{i=1}^r(1-\L^{1-d_i})^{-1}\\
&=|\pi_1(G)|
\L^{-\dim G}
  \prod_{i=1}^r(1-\L^{-d_i})^{-1}\prod_{i=1}^r(1-\L^{1-d_i})^{-1}\,, 
\end{align*}
by Proposition~\ref{groupprop}.
In view of the fact that $$Z(\pp^1,u)=(1-u)^{-1}(1-\L u)^{-1}\,,$$ this is 
Conjecture~\ref{general} for $C=\pp^1$.

\eject

\begin{appendix}

\section{The Motive of a Torsor}

\begin{center}
\small    By AJNEET DHILLON
\end{center}

\normalsize
\subsection{A Review of Voevodsky's Category of Motives}

We begin by briefly recalling the construction of the
triangulated category of effective geometrical motives
from \cite{voevodsky:00}. Denote by $\cor(k)$ the category
whose objects are schemes smooth over $k$ and
a morphism from $X$ to $Y$ is an algebraic cycle $Z$ on
$X\times Y$ such that each component of $Z$ is finite over $X$.
The finiteness condition allows one to define composition
without having to impose an adequate equivalence relation.
Note that $\cor(k)$ is an additive category with direct sum
given by $X\oplus Y = X \amalg Y$.

The homotopy category of bounded of bounded complexes
\[{\mathcal H}^b(\cor(k))\] is a triangulated category.
Let $T$ be the minimal thick subcategory containing
all complexes of the following two forms:

\noindent
1) $X\times{\mathbb A}^1 \stackrel{{\rm proj}}{\rightarrow} X$

\noindent
2) For every open cover $U,V$ of $X$ the complex
\[
U\cap V \hookrightarrow U \oplus V
\stackrel{(i_u,-i_v)}{\rightarrow} X.
\]

The triangulated category $\dmgm[\Z]$ of effective geometrical
motives is defined to be the Karoubian hull of the
localization of \[{\mathcal H}^b(\cor(k))\] with respect to $T$.
We will mostly be interested in its $\Q$-linearization
$\dmgm$. The obvious functor $\sm(k)\rightarrow\dmgm$ is denoted
$M_{gm}$. We now recall Voevodsky's alternative construction
of it.

A presheaf with transfers is a contravariant functor
on $\cor(k)$. It is called a sheaf with transfers if it
is a sheaf when restricted to the big etale site on
$\sm(k)$. We denote by ${\rm Shv}(\cor(k))$ the category of such sheaves.

 A presheaf with transfers $F$ is called
homotopy invariant if for all smooth schemes $X$,
the natural map $F(X)\rightarrow F(X\times\A^1)$
is an isomorphism. We denote by $\dmm[\Z]$ the
full subcategory of the derived category
$D^{-}({\rm Shv}(\cor(k)))$ consisting of those
complexes with homotopy invariant cohomology sheaves.
We will be mostly interested in its $\Q$-linearization
$\dmm$.

We denote by $\Delta^\bullet$ the cosimplicial
scheme with
\[
\Delta^n = \spec(k[x_0,x_1,\ldots,x_n]/\sum x_i =1)
\]
and face maps given by setting $x_i=0$.
Given a sheaf with transfers $F$ we denote by
$\cs(F)$ the complex associated to the simplicial
sheaf with transfers whose $n$th term is
\[
{\rm C}_n(F)(X) = F(X\times \Delta^n).
\]
Recall \cite[Lemma 3.2.1]{voevodsky:00} that the cohomology
sheaves of $\cs(F)$ are homotopy invariant.

\begin{theorem}
  The functor $\cs(-)$ extends to a functor
\[
{\bf R}\cs : D^-({\rm Shv}(\cor(k))) \rightarrow \dmm.
\]
This functor is left adjoint to the natural inclusion.
\end{theorem}

\begin{proof}
  See \cite[Theorem 3.2.3,\S 3.3]{voevodsky:00}.
\end{proof}

\begin{theorem}
For a perfect field $k$  there is a commutative
diagram of functors
\[
\xymatrix{
{\mathcal H}^b(\cor(k))\otimes\Q \ar[r]^(.4){L} \ar[d]
& D^-({\rm Shv}(\cor(k)))\otimes\Q \ar[d] \\
\dmgm \ar[r]^i & \dmm
}
\]
such that $i$ is a fully faithful embedding with dense image.
\end{theorem}

\begin{proof}
  See \cite[Theorem 3.2.6]{voevodsky:00} including the
construction of $i$.
\end{proof}

If the field $k$ admits a resolution of singularities
then there is a functor called the motive with
compact support:
\[
M_{gm}^c : {\rm sch}^{prop}/k \rightarrow \dmgm,
\]
here ${\rm sch}^{prop}/k$ is the category whose objects
are schemes of finite type over $k$ and morphisms
are proper maps.
If $Z$ is a closed subscheme of $X$ then there is
an exact triangle
\[
M_{gm}^c(Z)\rightarrow M_{gm}^c(X) \rightarrow
M_{gm}^c(X-Z)\rightarrow M_{gm}^c(Z)[1].
\]
For further properties see \cite[pg. 195]{voevodsky:00}.
We will now briefly recall the construction of
$M_{gm}^c$. For a scheme $X$ and a smooth scheme
$U$ define
$L^c(X)(U)$ to be the free abelian group
generated by closed integral subschemes of
$X\times U$ quasi-finite over $U$ and
dominant over a component of $U$. In this way we
obtain a sheaf
\[
L^c(X) :\cor(k)^{op}\rightarrow {\bf Ab}.
\]
We have a functor
\[
L^c:{\rm sch}^{prop}/k \rightarrow {\rm Shv}(\cor(k)).
\]
The motive with compact supports of $X$ is defined to be
${\bf R}\cs L^c(X)$. It is a theorem, see
\cite[Corollary 4.1.4]{voevodsky:00} that this sheaf
belongs to $\dmgm$.

\begin{proposition}
  Let $\Gamma$ be a finite group acting on the scheme
$X$ with quotient $Y=X/\Gamma$. Then
\[
L^c(X) \rightarrow L^c(Y)
\]
is a quotient in ${\rm Shv}(\cor(k))\otimes \Q$.
\end{proposition}

\begin{proof}
  The category ${\rm Shv}(\cor(k))\otimes \Q$ is equivalent
to the category of sheaves of $\Q$-vector spaces with
the equivalence being given by the functor $-\otimes\Q$.
For every $U$ the natural push forward map induces an
isomorphism
\[
L^c(X)(U)^\Gamma\rightarrow L^c(Y)(U),
\]
see \cite[1.7.6]{fulton:98}. The result
now follows.
\end{proof}

\begin{corollary}
  In the above notation the natural map
\[
L^c(X)\rightarrow L^c(Y)
\]
is a quotient in $D^-({\rm Shv}(\cor(k)))\otimes \Q$.
\end{corollary}

\begin{proof}
  Note that one can calculate maps from $L^c(X)$
to ${\mathcal F}$ in the derived category by taking
an injective resolution $I_\bullet$ of ${\mathcal F}$
and calculating homotopy classes of maps to $I_\bullet$.
So suppose
\[\phi:L^c(X)\rightarrow I_\bullet\]
is $\Gamma$-equivariant. So for each $\gamma\in \Gamma$
there exists
\[
h_\gamma: L^c(X)\rightarrow I_1
\]
with $\phi\circ \gamma -\phi = d_1 \circ h_\gamma$.
Then $\frac{1}{|\Gamma|}\sum_{\gamma\in\Gamma} \phi\circ\gamma$
is $\Gamma$ equivariant in ${\rm Shv}(\cor(k))\otimes \Q$
and this same map equals $\phi$ in the homotopy category.
Now apply the proposition.
\end{proof}

\begin{corollary}
\label{C:quot}
  In the same notation
\[
M^c_{gm}(X)\rightarrow M^c_{gm}(Y)
\]
is a quotient in $\dmgm$.
\end{corollary}

\begin{proof}
  This is because ${\bf R}\cs$ is  left adjoint to the
inclusion.
\end{proof}

\begin{proposition}
\label{P:a1same}
  Consider the family of inclusions
\[
i_s : X\hookrightarrow X\times{\mathbb A}^1.
\]
Then $\cs L^c(i_1) = \cs L^c(i_0)$ in
${\mathcal H}^b({\rm Shv}(\cor(k))\otimes \Q)$.
\end{proposition}

\begin{proof}
  This is well known. A proof can be found in
\cite[Lemma 2.17]{weibel:notes}.
\end{proof}

\subsection{The Main Result}

We assume throughout this section the characteristic
of the ground field $k$ is $0$.

\begin{theorem}
\label{T:deform}
  Let $G$ be a split semisimple connected group or a 
connected unipotent
group over $k$. Let $x\in G(A)$ where $A$
is a finite generated $k$ algebra that is a domain.
Then there is an open affine $\spec(A')\subseteq\spec(A)$
and a finite Galois cover $\spec(B)\rightarrow\spec(A')$
and $y\in G(B[t])$ such that

\noindent
(i) $y(0)$ is the constant morphism to the identity

\noindent
(ii) The following diagram commutes
\[
\xymatrix{
\spec(A') \ar[r]  & \spec A \ar[r]^x & G \\
\spec B \ar[urr]_{y(1)} \ar[u]&
}
\]
\end{theorem}

\begin{proof}
  The result is straightforward in the case where
$G$ is a connected unipotent group as in this
case the underlying variety of $G$ is an
${\mathbb A}^n$, see \cite[pg. 243]{springer:98}.

So we assume $G$ is semisimple.
 Let $\tilde{G}$ be the
universal cover of $G$. As
\[
\tilde{G}\rightarrow G
\]
is Galois, by replacing $\spec(A)$ by
$\spec(A)\times_{G,x} \tilde{G}$ we may assume
that $G$ is simply connected.

According to
\cite{steinberg:62} there is a unipotent group $U$
and a morphism
\[
\phi:U\rightarrow G
\]
that is surjective on $L$ points for every field $L$.
Let $K$ be the function field of $A$ and $x_K$ be the
$A$-point of $G$ restricted to $K$. There is a
$K$-point $x'$ of $U$ mapping to $x_K$ via $\phi$.
By examining denominators we can find an open affine
\[
\spec(A')\hookrightarrow\spec(A)
\]
such that $x_{A'}$ lifts to a $A'$-point $z$ of $U$.
Now the underlying variety of $U$ is again
an ${\mathbb A}^n$. 
\end{proof}

Let $\pi:Y\rightarrow X$ be a finite Galois cover
with Galois group $\Gamma$. Let $P$ be a $G$-torsor
trivialized by $Y$, that is $Y\times_X P\isoarrow Y\times G$.
The action of $\Gamma$ on $Y$ lifts to an action
of $\Gamma$ on $Y\times G$ with quotient $P$. This
action is determined by a 1-cocyle
\[
n:\Gamma\times Y \rightarrow G.
\]

\begin{proposition}
\label{P:deftriv}
  Let $G$ be a connected split semisimple group or a group
whose underlying variety is ${\mathbb A}^n$.
 Let $P\rightarrow X$ be a $G$-torsor.
Then there is an open affine $X'\subseteq X$ and a Galois
cover $Y'\rightarrow X'$ trivializing $P$. Furthermore
if $\Gamma$ is the Galois group and the cocycle
\[
n:\Gamma \times Y \rightarrow G
\]
defines the action we may assume that
$n$ extends to
\[
n_t:\Gamma \times Y \times {\mathbb A}^1 \rightarrow G
\]
with $n_0 = n$ and $n_1$ constant at the identity.
\end{proposition}

\begin{proof}
  The first part is standard, using Zariski's main theorem.
The second part is by repeated applications of \ref{T:deform}.
\end{proof}

We denote by $K_0(\dmgm)$ the $K$-group of the triangulated
category $\dmgm$. It is the free abelian group
on the objects of $\dmgm$ subject to the relations
\[
Y = X + Z \quad\text{for each exact triangle}\quad
X\rightarrow Y\rightarrow Z\rightarrow X[1].
\]
The tensor product of the category $\dmgm$
makes $K_0(\dmgm)$ into a ring.
We have a ring homomorphism
\[
\chi^c_{\rm Mot}:K_0(Var_k)\rightarrow K_0(\dmgm)
\]
given by the motive with compact supports.
We denote by $\eumc(X)$ the image of the variety
$X$ under this homomorphism.

\begin{theorem}
  Let $G$ be a connected linear algebraic group
that is split over $k$. Let $P$ be a $G$-torsor
over $X$ with $X$ of finite type.
Then $\eumc(P)=\eumc(X)\eumc(G)$.
\end{theorem}

\begin{proof}
  By noetherian induction it suffices to find
an open subset $X'$ of $X$ such that
\[
\eumc(X')\eumc(G) = \eumc(P|_U).
\]
First we assume that $G$ is as in \ref{P:deftriv}.
Then we can find $\Gamma$, $n_t$ $X'$ and  $Y'$ as in the proposition.
Now the natural map
\[
M^c_gm(P\times_X Y')\rightarrow  M^c_{gm}(P|_{X'})
\]
is a quotient by \ref{C:quot}. On the other
hand by \ref{P:a1same} the cocycle is
$M^c_{gm}(n)$ is trivial. Hence the result.

For a general group note that $P\rightarrow P/{R_u(G)}$
is a $R_u(G)$-torsor.  So we may assume $G$
is split  reductive. But then $R(G)$ is a torus and
every torsor for a torus is Zariski trivial so we
are reduced to the semisimple case.
\end{proof}

\end{appendix}

\bibliographystyle{alpha}
\bibliography{./../../latestbib/alggrp.bib,./../../latestbib/mybib.bib,./../../latestbib/sga.bib,./../../latestbib/mybib2.bib}

\end{document}